\documentclass[graybox]{svmult}

\usepackage{mathptmx}       
\usepackage{helvet}         
\usepackage{courier}        
\usepackage{type1cm}        
\usepackage{makeidx}         
\usepackage{graphicx}        
\usepackage{multicol}        
\usepackage[bottom]{footmisc}

\usepackage[T1]{fontenc}
\usepackage{graphicx}
\usepackage{amsmath,amssymb,mathtools}
\usepackage{bm}
\usepackage{booktabs}
\usepackage{enumitem}
\usepackage{algorithm}
\usepackage{algpseudocode}
\usepackage{pgfplots}
\usepackage{pgfplotstable}
\usepackage{longtable}
\usepackage{booktabs}
\usepackage{longtable}

\pgfplotsset{compat=1.18}

\newcommand{\Bset}{\mathcal{B}}
\newcommand{\X}{\mathcal{X}}
\newcommand{\R}{\mathbb{R}}
\newcommand{\eps}{\varepsilon}

\makeindex             

\begin{document}

\title*{Global Optimization of Mixed-Integer Bound-Constrained Quadratic Fractional Problems}
\titlerunning{Global Optimization of Fractional MIQPs}
\author{Marco Locatelli \and Veronica Piccialli \and Antonio M. Sudoso}
\institute{Marco Locatelli \at Dipartimento di Ingegneria e Architettura, Universit\`{a} di Parma, Parco Area delle Scienze, 181/A,
Parma 43124, Italy \email{marco.locatelli@unipr.it}
\and Veronica Piccialli \at Dipartimento di Ingegneria Informatica, Automatica e Gestionale "A. Ruberti", Sapienza Universit\`{a}
di Roma, Via Ariosto 25, Rome 00185, Italy \email{veronica.piccialli@uniroma1.it} \and Antonio M. Sudoso \at Dipartimento di Matematica "Tullio Levi-Civita", Universit\`{a} degli Studi di Padova, Via Trieste 63, Padua 35121, Italy \email{antoniomaria.sudoso@unipd.it}}
%
%
\maketitle

\abstract*{We study mixed-integer quadratic fractional optimization problems with binary and continuous box-constrained variables. The objective is the ratio of two quadratic functions, a structure that arises naturally in portfolio selection, signal processing, and data analytics. We propose an exact algorithm based on Dinkelbach's parametric reformulation. At each Dinkelbach iteration, the fractional problem is replaced by a nonconvex mixed-integer quadratic subproblem with box constraints. This subproblem is solved globally by a specialized branch-and-bound solver. To accelerate the outer scheme, we use a Variable Neighborhood Search heuristic to provide an initial feasible solution to the Dinkelbach scheme. The heuristic alternates between continuous local searches and a binary-space exploration phase in which the continuous variables of a feasible incumbent are fixed. Numerical results show that the proposed approach outperforms general-purpose exact algorithms based on Dinkelbach's parametric reformulation.}

\abstract{We study mixed-integer quadratic fractional optimization problems with binary and continuous box-constrained variables. The objective is the ratio of two quadratic functions, a structure that arises naturally in portfolio selection, signal processing, and data analytics. We propose an exact algorithm based on Dinkelbach's parametric reformulation. At each Dinkelbach iteration, the fractional problem is replaced by a nonconvex mixed-integer quadratic subproblem with box constraints. This subproblem is solved globally by a specialized branch-and-bound solver. To accelerate the outer scheme, we use a Variable Neighborhood Search heuristic to provide an initial feasible solution to the Dinkelbach scheme. The heuristic alternates between continuous local searches and a binary-space exploration phase in which the continuous variables of a feasible incumbent are fixed. Numerical results show that the proposed approach outperforms general-purpose exact algorithms based on Dinkelbach's parametric reformulation.}

\section{Introduction}
Many optimization models evaluate a feasible solution through a normalized performance measure rather than an absolute objective value. Such criteria arise when the quantity of interest is an efficiency, density, average cost, or tradeoff between two competing terms. In these settings, the objective typically balances a numerator measuring loss, disagreement, risk, or resource consumption against a denominator measuring scale, volume, benefit, or signal strength.

This structure appears in a wide range of applications. Within this literature, quadratic fractional programming occupies a particularly important niche because the numerator and denominator can both model energy, risk, size, density, or normalized utility. 
In clustering and community detection, the quality of a partition is often measured by a cut value normalized by volume or cardinality \cite{piccialli2023global}, leading to criteria such as ratio cut and normalized cut \cite{chan1993spectral,de2022ratio}. In biclustering and data mining, submatrices may be assessed through density, homogeneity, or normalized disagreement \cite{fan2012multi,sudoso2025semidefinite,sudoso2026exact}. In portfolio optimization, risk-adjusted performance measures such as the Sharpe ratio compare return with risk \cite{sharpe1966mutual}. Ratio objectives also arise in signal processing and wireless communications, where signal-to-interference-plus-noise ratios measure transmission quality \cite{cai2014approach}. Mixed-integer optimization models involving ratios of linear functions also arise in the context of cyclic process scheduling problems \cite{you2009dinkelbach}.

These examples lead naturally to fractional optimization problems where the numerator and the positive denominator are nonlinear functions. {Indeed, in all these applications, positivity of the denominator arises naturally. In normalized-cut and related graph-partitioning models, the denominator represents the cardinality or volume of a cluster and is therefore positive for nonempty clusters. In portfolio models, the denominator typically represents a risk measure, which is positive when zero-risk portfolios are excluded. In wireless communication models based on signal-to-noise or signal-to-interference-plus-noise ratios, the denominator contains noise power and is therefore positive.}
Despite their modeling value, fractional objectives are analytically and computationally challenging. Even simple choices of numerator and denominator can produce nonconvex problems, and the ratio may fail to inherit useful properties of its components. 
Specialized approaches have been defined for fractional linear problems \cite{benson2007simplicial,benson2004global,jiao2015practicable,luo2025new}.

In this paper, we consider the mixed-integer bound-constrained fractional quadratic program (Fractional MIQP) of the form
\begin{equation}\label{eq:MIQP}
    \min_{x \in \mathcal{X}} \frac{N(x)}{D(x)},
\end{equation}
with
\[
N(x):=x^\top A x+a^\top x+a_0,\quad D(x):=x^\top B x+b^\top x+b_0
\]
and $D(x)>0$ for all $x\in\mathcal X$ and
\begin{equation*}
\label{eq:X}
\X=\left\{x\in\R^n:0\le x_i\le1\ \forall i\notin\Bset,\ \ x_i\in\{0,1\}\ \forall i\in\Bset\right\},
\end{equation*}
where $\Bset\subseteq\{1,\dots,n\}$ is the binary index set. 
We assume $A \in \mathbb{S}^n$ and $B \in \mathbb{S}^n$ indefinite, $D(x)>0$ on $\X$ and compactness of $\X$, which is immediate under box and binary constraints.  This model is challenging due to the simultaneous presence of a nonconvex objective function and binary variables coupled with continuous box-constrained variables.
 
Classical fractional programming methods address these difficulties through parametric reformulations, as in Dinkelbach-type methods \cite{dinkelbach1967nonlinear,schaible1976fractional}. The approach solves a sequence of subproblems, that remain difficult when the numerator, denominator, and the feasible set are nonconvex.
Besides the Dinkelbach-type method, there are also various methods, including the Charnes–Cooper transformation \cite{charnes1962programming}, and proximal-gradient algorithms \cite{boct2017proximal}. 
The proximal gradient method proposed in \cite{boct2017proximal} can be applied in the continuous case, when the numerator function is proper, convex, and lower semicontinuous, and the denominator is differentiable with Lipschitz continuous gradient and either concave or convex.
The Charnes-Cooper transformation \cite{charnes1962programming} transforms the problem into a parameter-free equivalent problem, which is in general nonconvex, apart from the special case where $N(x)$ is convex and $D(x)$ is concave.
Furthermore, convexification techniques can be used both in the continuous and in the discrete setting \cite{benson2002using,bugarin2016minimizing,he2025convexification,de2026beyond} to design exact approaches. 

 For fractional linear 0--1 and mixed-integer linear fractional problems, conic and MILP reformulations have been explored, especially when structural properties allow tighter bounds \cite{mehmanchi2019fractional,you2009dinkelbach}. 
In the continuous quadratic case, various branch-and-bound and branch-and-reduce approaches have been proposed \cite{qu2007efficient,shen2009solving}.  On the heuristic side, local-search and variable-neighborhood methods remain attractive because they can exploit problem structure and rapidly improve primal incumbents, particularly when incremental updates make neighborhood evaluation cheap \cite{consolini2020efficient,mladenovic1997variable}.
 Recently, an exact method for ternary fractional QPs has been proposed in \cite{de2026beyond}, exploiting semidefinite programming and the Dinkelbach algorithm.  To the best of our knowledge, no exact approach has been proposed so far for mixed-integer quadratic fractional optimization problems with bound constraints and both binary and continuous variables. This gap motivates the development of a tailored global optimization method for this class of problems.

The main contributions of the paper are the following.
\begin{itemize}
    \item We propose an exact solution method for fractional MIQPs based on
    Dinkelbach's parametric reformulation. At each iteration, the original
    fractional problem is replaced by a nonconvex mixed-integer quadratic
    subproblem, which is solved to global optimality. We show that these
    subproblems can be reformulated as box-constrained quadratic programs
    and solved efficiently by a tailored approach.

    \item We introduce a hybrid heuristic that combines
    continuous local optimization with Variable Neighborhood Search on the
    binary variables. The heuristic provides high-quality feasible solutions
    that are used to initialize the Dinkelbach scheme, thereby reducing the
    number of subproblems to be solved to global optimality.

    \item We provide a computational study on benchmark instances varying in
    problem size, matrix density, and percentage of binary variables. The
    results show that the proposed approach globally solves all tested
    fractional MIQP instances and outperforms a
    general-purpose Dinkelbach implementation based on Gurobi, especially on
    larger and denser instances.
\end{itemize}

The remainder of the paper is organized as follows. Section~\ref{sec:method}
presents the Dinkelbach reformulation, the BoxQP reformulation of the
parametric subproblems, and the hybrid initialization heuristic.
Section~\ref{sec:results} describes the benchmark instances and reports the
computational results. Finally, Section~\ref{sec:concl} concludes the paper
and discusses directions for future work.

\section{Methodology}\label{sec:method}

Our solution approach is based on Dinkelbach's algorithm. Dinkelbach introduced a classical framework for fractional programming~\cite{dinkelbach1967nonlinear}, which was later extended to more general settings~\cite{rodenas1999extensions}. The key idea is to replace the direct optimization of a ratio by the iterative solution of a sequence of parametric optimization problems.
We define the parametric value function
\begin{equation}
\label{eq:F}
F(\lambda)
=
\min_{x\in\X}
\left\{
N(x)-\lambda D(x)
\right\},
\qquad \lambda\in\R .
\end{equation}
Since $F$ is the pointwise minimum of affine functions of $\lambda$, it is concave. Moreover, under the standing assumption that $D(x)>0$ for all $x\in\X$, the function $F$ is strictly decreasing. For any feasible solution $x\in\X$ and any value $\lambda=N(x)/D(x)$, we have
\[
F(\lambda)
\le
N(x)-\lambda D(x)
=
0.
\]
Therefore, the optimal value of the fractional problem is characterized as the unique zero of $F$:
\[
\lambda^\star
=
\min_{x\in\X}\frac{N(x)}{D(x)}
\quad
\Longleftrightarrow
\quad
F(\lambda^\star)=0.
\]
Thus, instead of solving the fractional problem directly, Dinkelbach's method iteratively evaluates $F(\lambda)$ and updates the parameter $\lambda$ using the ratio value associated with the optimizer of the current subproblem. Note that the optimizer of each parametric subproblem must be computed globally. If only a local solution were used, the computed value would not necessarily coincide with $F(\lambda)$, and the zero-finding interpretation of Dinkelbach's method would no longer provide a certificate of global optimality.

In our fractional MIQP setting, each Dinkelbach iteration requires solving, for a fixed value of $\lambda$, the parametric subproblem
\begin{equation}
\label{eq:inner_qp}
    \min_{x \in \X}
    x^\top Q(\lambda)x + q(\lambda)^\top x + q_0(\lambda),
\end{equation}
where
\[
    Q(\lambda)=A-\lambda B,\qquad
    q(\lambda)=a-\lambda b,\qquad
    q_0(\lambda)=a_0-\lambda b_0.
\]
Problem~\eqref{eq:inner_qp} is a nonconvex MIQP. Hence, at each iteration, our method solves a nonconvex MIQP to global optimality. The globally optimal solution of this subproblem is then used to update the parameter $\lambda$ and to evaluate the Dinkelbach residual. The procedure terminates when this residual is below the prescribed tolerance. 
The appeal of Dinkelbach's reformulation is that it avoids handling the fractional objective explicitly. Instead, the original problem is transformed into a sequence of optimization problems with a quadratic objective function. 
However, in the present setting, these subproblems remain computationally challenging because they are nonconvex MIQPs. Consequently, the efficiency of the overall method strongly depends on the ability to solve the parametric subproblems globally and efficiently. This motivates the use of a tailored procedure.
The fix-and-bound method proposed in~\cite{locatelli2025fix} is designed for large-scale nonconvex box-constrained quadratic programming (BoxQP) problems. Its main idea is to combine a branch-and-cut framework with strong lower bounds based on semidefinite programming and variable-fixing mechanisms. We exploit this solver to handle the parametric subproblems arising within Dinkelbach's method.
To this end, we reformulate the mixed-integer structure of Problem~\eqref{eq:inner_qp} so that it can be treated within a box-constrained quadratic programming framework. For BoxQP problems, if the diagonal coefficient associated with variable $x_i$ satisfies $Q(\lambda)_{ii}\le 0$ then, by the optimality conditions for BoxQP, there exists a global minimizer in which $x_i$ is binary, even if $x_i$ is modeled as a continuous variable in $[0,1]$~\cite{hansen1993global}. Conversely, for each binary variable \(x_j\), \(j\in\mathcal B\), its
binary nature can be enforced through a diagonal-shifting reformulation. Let $\sigma_j \ge \max\{Q(\lambda)_{jj},0\}$,
and define
\[
\widetilde Q(\lambda)_{jj}
=
Q(\lambda)_{jj}-\sigma_j,
\qquad
\widetilde q(\lambda)_j
=
q(\lambda)_j+\sigma_j,
\qquad
\widetilde q_0(\lambda)=q_0(\lambda).
\]
All other entries of \(Q(\lambda)\) and \(q(\lambda)\) are left unchanged.
Then, for any \(x_j\in\{0,1\}\),
\[
\widetilde Q(\lambda)_{jj}x_j^2
+
\widetilde q(\lambda)_j x_j
=
Q(\lambda)_{jj}x_j^2
+
q(\lambda)_j x_j
+
\sigma_j(x_j-x_j^2)
=
Q(\lambda)_{jj}x_j^2
+
q(\lambda)_j x_j,
\]
because \(x_j^2=x_j\). Hence the reformulation preserves the objective
value on the original mixed-integer feasible set, while ensuring that
\(\widetilde Q(\lambda)_{jj}\le 0\). Therefore, by the BoxQP property used
in~\cite{hansen1993global,locatelli2025fix}, the transformed box-constrained
quadratic problem admits a global minimizer in which the variables
\(x_j\), \(j\in\mathcal B\), are binary.
As a result, each parametric MIQP subproblem generated by Dinkelbach's algorithm is transformed into a BoxQP and solved globally by the fix-and-bound solver. The overall algorithm therefore combines an outer Dinkelbach scheme, which updates the fractional parameter, with an inner global optimization procedure, which solves each parametric quadratic subproblem to global optimality.

Schaible~\cite{schaible1976fractional} proved the superlinear convergence of Dinkelbach's algorithm. A similar convergence result can be established for the variant of Dinkelbach's algorithm applied to fractional MIQPs. The pseudocode of the method is given in Algorithm~\ref{pc:dink}. The parameter $\eps$ denotes the stopping tolerance: the algorithm terminates once the Dinkelbach residual is at most $\eps$, thereby returning a solution of Problem~\eqref{eq:MIQP} within the prescribed optimality tolerance.

\begin{algorithm}[!ht]
\caption{Exact algorithm for fractional MIQP}\label{pc:dink}
\begin{algorithmic}[1]
\Require Initial point $x^0 \in \mathcal{X}$, tolerance $\eps > 0$
\State Set $\lambda^0 \gets N(x^0)/D(x^0)$
\For{$t=0,1,\dots$}
  \State Build $(\widetilde Q(\lambda^t),\widetilde q(\lambda^t),
  \widetilde q_0(\lambda^t))$ and globally solve the corresponding
  BoxQP reformulation of \eqref{eq:inner_qp} by fix-and-bound, obtaining
  \[
      x^{t+1} \in 
      \arg\min_{x \in \mathcal{X}}
      \left\{ N(x) - \lambda^t D(x) \right\}.
  \]
  \State $r^t \gets N(x^{t+1})-\lambda^tD(x^{t+1})$
  \If{$|r^t|\le \eps$}
      \State \Return $x^{t+1},\, \lambda^t$
  \EndIf
  \State $\lambda^{t+1}\gets N(x^{t+1})/D(x^{t+1})$
\EndFor
\end{algorithmic}
\end{algorithm}

Since solving a single parametric subproblem \eqref{eq:inner_qp} to global
optimality is computationally expensive, obtaining a strong initial feasible
solution is important in practice. A high-quality incumbent provides a good
initial value of the Dinkelbach parameter \(\lambda\), and may significantly
reduce the number of MIQP subproblems that need to be solved. In the ideal
case where the heuristic solution is already globally optimal, only one
parametric MIQP subproblem is required to certify optimality.
For this reason, we design a primal heuristic to generate high-quality feasible
solutions before starting the exact Dinkelbach scheme. The heuristic is based on
Variable Neighborhood Search (VNS), a well-established metaheuristic for
combinatorial and mixed-integer optimization introduced in
\cite{hansen2010variable,mladenovic1997variable}. The main idea of VNS is to systematically explore
increasingly distant neighborhoods of the incumbent solution in order to escape
local minima and improve diversification. The method alternates between a
shaking phase, which perturbs the current solution by random moves, and a local
search phase, which intensifies the search around the perturbed solution. In the
context of Problem~\eqref{eq:MIQP}, we combine VNS on the binary variables with
continuous local optimization over the box-constrained variables. Starting from a random point in the box $[0,1]^n$, the method first applies a continuous nonlinear local solver to the relaxation of the problem, producing a KKT point. The binary variables are then rounded to obtain a feasible mixed-integer solution.
After initialization, the algorithm alternates between two phases. In the first phase, a VNS procedure explores the neighborhood of the current solution by modifying only the binary variables, while all continuous variables remain fixed. The neighborhood structure is based on binary flips and is explored through the classical shaking and local-search strategy of VNS. In the second phase, the binary variables are fixed and a continuous local optimization problem is solved over the remaining continuous variables. This alternating strategy allows the heuristic to jointly exploit combinatorial exploration on the discrete variables and continuous nonlinear refinement on the box-constrained variables.
The overall procedure is embedded in a multi-start framework, where several random initial points are independently refined. The complete heuristic is summarized in Algorithm~\ref{alg:hybrid_heuristic}.

\begin{algorithm}[!ht]
\caption{Hybrid local search heuristic for fractional MIQP}
\label{alg:hybrid_heuristic}
\begin{algorithmic}[1]
\Require Initial point $x^0 \in [0,1]^n$
\State Find a KKT point of $N(x)/D(x)$ over $[0,1]^n$ starting from $x^0$ using continuous local search
\State Round $x_i$ to nearest integer in $\{0,1\}$ for all $i \in \Bset$.
\State $x^\star \gets x$ (current point after rounding)
\State Initialize $N(x^\star), D(x^\star)$
\State improved $\gets 1$
\While{improved $= 1$}
  \State improved $\gets 0$
  \State \textbf{Phase 1: Binary variable neighborhood search}
  \State Fix $x_j = x^\star_j$ for all $j\not \in \Bset$
  \State $x \gets \mathrm{VNS}(x^\star)$
  \If{$N(x)/D(x) < N(x^\star)/D(x^\star)$}
    \State $x^\star \gets x$
    \State improved $\gets 1$
  \EndIf
  \State \textbf{Phase 2: Continuous variable optimization}
    \State Fix $x_j = x^\star_j$ for all $j \in \Bset$
    \State Find a KKT point of $N(x)/D(x)$ optimizing over $x_i \in [0,1]$ for all $i \notin \Bset$ by a continuous local search starting from $x^\star$
    \If{$N(x)/D(x) < N(x^\star)/D(x^\star)$}
      \State $x^\star \gets x$
      \State improved $\gets 1$
    \EndIf
\EndWhile
\State \Return $x^\star$
\end{algorithmic}
\end{algorithm}

\section{Results}\label{sec:results}
In this section, we evaluate the computational performance of the proposed solution approach on a set of randomly generated instances. We first describe the instance-generation procedure and the implementation details. We then compare the proposed \texttt{F\&B-Ratio} method with Gurobi in terms of solution quality, runtime, and robustness across the test set.
The \texttt{F\&B-Ratio} method is implemented in C++. All experiments are performed on a macOS system equipped with an Apple M4 Max chip with 14 cores and 36 GB of RAM, running macOS 15.5.

\subsection{Instance generation}\label{sec:inst}

To evaluate the computational performance of the proposed approach, we generate a collection of random fractional MIQP instances with controlled problem size, sparsity, and integrality structure. Each instance is characterized by the dimension $n$, the binary-variable percentage $p$, the density level $d$, and a random seed $s$.
For each test instance, the numerator coefficients $(A,a,a_0)$ and denominator coefficients $(B,b,b_0)$ are generated independently at random. The quadratic matrices $A,B\in\mathbb{S}^n$ are sparse. Given a target density $d\in[0,100]$, each entry in the lower triangular part is set to a nonzero value with probability $d/100$ and to zero otherwise. Symmetry is then enforced by mirroring the lower triangular part to the upper triangular part. Whenever an entry is nonzero, it is sampled uniformly from the integer interval $\{-50,\ldots,50\}$. The linear vectors $a,b\in\mathbb{R}^n$ and the constant term $a_0$ are generated independently with i.i.d. entries uniformly distributed over $\{-50,\ldots,50\}
$. The denominator constant is set to
\begin{equation*}
    b_0
=
1
+
\sum_{i=1}^n |B_{ii}|
+
2\sum_{1\le i<j\le n} |B_{ij}|
+
\sum_{i=1}^n |b_i|,
\end{equation*}
which ensures that the denominator satisfies $D(x)\ge 1$ over the feasible domain.
To introduce mixed-integer structure, a subset $\mathcal{B}\subseteq\{1,\ldots,n\}$ of binary variables is selected uniformly at random without replacement. For a prescribed binary percentage $p$, we set $|\mathcal{B}| = \left\lfloor p n \right\rfloor,
$ and impose $x_i\in\{0,1\}$ for $i\in\mathcal{B}$ and $x_i\in[0,1]$ for $i\notin\mathcal{B}$. Thus, the remaining variables are continuous and bounded in $[0, 1]$. A summary of the generated instances is reported in Table~\ref{tab:grid}. The instances are publicly available at \url{https://github.com/antoniosudoso/fractional-MIQP}

\begin{table}[!ht]
\centering
\caption{Experimental design for the generated MIQP-Ratio benchmark instances. The design varies the problem dimension, matrix density, percentage of binary variables, and random seed, resulting in 108 instances in total.}
\label{tab:grid}
\begin{tabular}{lcc}
\toprule
Component & Values & Count \\
\midrule
Problem size $n$ & $\{70,80,90,100,110,120\}$ & 6 \\
Density percentage $d$ & $\{25\%,50\%,75\%\}$ & 3 \\
Binary-variable percentage $\rho$ & $\{25\%,50\%\}$ & 2 \\
Seeds per setting & $\{1,2,3\}$ & 3 \\
\midrule
Total instances &  & $6\times3\times2\times3=108$ \\
\bottomrule
\end{tabular}
\end{table}

\subsection{Computational results}

Before starting the Dinkelbach procedure, we compute an initial solution
$x^0\in\mathcal{X}$ by running the heuristic in
Algorithm~\ref{alg:hybrid_heuristic} in a multi-start fashion from 100 randomly generated initial points. Among the solutions obtained, we select the one with the lowest fractional objective value and use the corresponding ratio as the initial value of the Dinkelbach parameter. As a continuous local solver, we use SNOPT version~7.7~\cite{gill2005snopt}. The stopping tolerance in Algorithm~\ref{pc:dink} is set to $\epsilon=10^{-6}$.

The proposed method, denoted by \texttt{F\&B-Ratio}, solves each parametric MIQP
subproblem arising in the Dinkelbach scheme by means of the fix-and-bound solver of~\cite{locatelli2025fix}, using its default parameter settings. As a benchmark, we consider the same Dinkelbach scheme in which each nonconvex MIQP
subproblem is solved by Gurobi version~12.0.2 \cite{gurobi} with default parameters. Thus, both methods use the same outer fractional programming framework and differ only in the solver used to handle the parametric quadratic subproblems. To the best of our knowledge, there is no specialized global optimization algorithm directly available for Problem~\eqref{eq:MIQP}; hence, this comparison provides a natural benchmark against a state-of-the-art commercial solver for nonconvex MIQPs.

All runs are performed with a time limit of 3600 seconds. Table~\ref{tab:gurobi-fnb-full} reports the detailed results for all benchmark instances. For each instance, we show the initial objective value obtained by the multi-start heuristic, the certified optimal objective value, the total runtime, and the number of
Dinkelbach iterations required by each method. The results show that \texttt{F\&B-Ratio} is substantially more robust on the full test
set. It solves all 108 instances within the time limit, whereas the Gurobi-based
implementation solves 82 instances and reaches the time limit on the remaining
26. The difference becomes more pronounced as the instances become larger and
denser. In particular, Gurobi solves all instances with $n=70$ and $n=80$, but
the number of time-limit runs increases for $n\geq 90$, especially for the
densest instances. By contrast, \texttt{F\&B-Ratio} remains effective across all tested
values of $n$, $p$, and $d$.

The number of Dinkelbach iterations is small for both implementations, ranging
between 1 and 3 over the entire benchmark set.  This behavior is consistent with
the quality of the initialization provided by the multi-start heuristic. In
particular, the heuristic already identifies a globally optimal solution on
16 instances. Over the full test set, the relative optimality gap of the initial solution compared with respect to the optimal solution, ranges from $0 \%$ to $7.34\%$, with an average gap of $1.91\%$ and a median gap of
$0.84\%$. Hence, the Dinkelbach procedure typically starts from a high-quality ratio value and requires only a few outer iterations to certify global optimality. {Notably, the computational time of the hybrid local search heuristic is negligible, taking only a few milliseconds.}
Finally, we note that, for instances solved to optimality, Gurobi and \texttt{F\&B-Ratio} perform the same number of iterations. By contrast, when Gurobi reaches the time limit, its iteration count may be lower; see, for example, the instance with $n=100$, $p=50$, $d=75$, and $s=2$.

{
\footnotesize
\setlength{\tabcolsep}{4pt}
\begin{longtable}{rrrrrrrrrr}
\caption{Comparison of Gurobi and \texttt{F\&B-Ratio} on all benchmark instances. For each instance, the table reports the initial objective value, the optimal objective value, total runtime in seconds, and number of Dinkelbach iterations.}
\label{tab:gurobi-fnb-full}\\

\toprule
\multicolumn{6}{c}{Instance} 
& \multicolumn{2}{c}{Gurobi} 
& \multicolumn{2}{c}{\texttt{F\&B-Ratio}} \\
\cmidrule(lr){1-6}
\cmidrule(lr){7-8}
\cmidrule(lr){9-10}
$n$ & $p$ & $d$ & $s$ & Init. obj. & Opt. obj.
& Time (s) & Iter.
& Time (s) & Iter. \\
\midrule
\endfirsthead

\toprule
\multicolumn{6}{c}{Instance} 
& \multicolumn{2}{c}{Gurobi} 
& \multicolumn{2}{c}{\texttt{F\&B-Ratio}} \\
\cmidrule(lr){1-6}
\cmidrule(lr){7-8}
\cmidrule(lr){9-10}
$n$ & $p$ & $d$ & $s$ & Init. obj. & Opt. obj.
& Time (s) & Iter.
& Time (s) & Iter. \\
\midrule
\endhead

\midrule
\multicolumn{10}{r}{Continued on next page}\\
\endfoot

\bottomrule
\endlastfoot

70 & 25 & 25 & 1 & -2.1005 & -2.1710 & 0.73 & 2 & 20.22 & 2 \\
70 & 25 & 25 & 2 & -1.3065 & -1.3154 & 2.53 & 2 & 49.96 & 2 \\
70 & 25 & 25 & 3 & -1.9588 & -1.9588 & 0.44 & 1 & 10.31 & 1 \\
70 & 25 & 50 & 1 & -1.1667 & -1.1928 & 27.19 & 3 & 76.17 & 3 \\
70 & 25 & 50 & 2 & -0.9003 & -0.9125 & 19.83 & 2 & 39.80 & 2 \\
70 & 25 & 50 & 3 & -1.1748 & -1.1955 & 10.54 & 2 & 35.23 & 2 \\
70 & 25 & 75 & 1 & -0.9756 & -0.9789 & 108.89 & 2 & 110.11 & 2 \\
70 & 25 & 75 & 2 & -0.6321 & -0.7056 & 718.74 & 2 & 96.37 & 2 \\
70 & 25 & 75 & 3 & -1.0203 & -1.0233 & 39.98 & 2 & 53.80 & 2 \\
70 & 50 & 25 & 1 & -2.1564 & -2.1710 & 0.73 & 2 & 19.06 & 2 \\
70 & 50 & 25 & 2 & -1.2877 & -1.3154 & 2.91 & 2 & 51.06 & 2 \\
70 & 50 & 25 & 3 & -1.9588 & -1.9588 & 0.47 & 1 & 11.00 & 1 \\
70 & 50 & 50 & 1 & -1.1712 & -1.1928 & 16.14 & 2 & 81.11 & 2 \\
70 & 50 & 50 & 2 & -0.8693 & -0.9125 & 17.88 & 2 & 38.39 & 2 \\
70 & 50 & 50 & 3 & -1.1905 & -1.1955 & 11.41 & 2 & 35.46 & 2 \\
70 & 50 & 75 & 1 & -0.9765 & -0.9789 & 86.48 & 2 & 77.32 & 2 \\
70 & 50 & 75 & 2 & -0.6851 & -0.7056 & 583.13 & 2 & 98.15 & 2 \\
70 & 50 & 75 & 3 & -1.0203 & -1.0233 & 36.90 & 2 & 46.80 & 2 \\
80 & 25 & 25 & 1 & -1.8984 & -1.9175 & 2.71 & 3 & 25.88 & 3 \\
80 & 25 & 25 & 2 & -1.4282 & -1.4283 & 3.18 & 2 & 55.06 & 2 \\
80 & 25 & 25 & 3 & -1.8923 & -1.8923 & 0.58 & 1 & 12.89 & 1 \\
80 & 25 & 50 & 1 & -1.1612 & -1.1612 & 18.79 & 1 & 35.07 & 1 \\
80 & 25 & 50 & 2 & -0.8626 & -0.8626 & 20.62 & 1 & 24.02 & 1 \\
80 & 25 & 50 & 3 & -1.1478 & -1.1586 & 28.22 & 2 & 56.89 & 2 \\
80 & 25 & 75 & 1 & -0.9207 & -0.9525 & 186.44 & 2 & 48.84 & 2 \\
80 & 25 & 75 & 2 & -0.6802 & -0.7065 & 1,275.71 & 2 & 54.32 & 2 \\
80 & 25 & 75 & 3 & -0.8620 & -0.8779 & 378.17 & 3 & 76.47 & 3 \\
80 & 50 & 25 & 1 & -1.9175 & -1.9175 & 0.99 & 1 & 12.60 & 1 \\
80 & 50 & 25 & 2 & -1.4151 & -1.4283 & 3.38 & 2 & 56.30 & 2 \\
80 & 50 & 25 & 3 & -1.8901 & -1.8923 & 1.23 & 2 & 25.06 & 2 \\
80 & 50 & 50 & 1 & -1.1547 & -1.1612 & 36.39 & 2 & 73.04 & 2 \\
80 & 50 & 50 & 2 & -0.8110 & -0.8626 & 48.17 & 2 & 47.98 & 2 \\
80 & 50 & 50 & 3 & -1.1277 & -1.1586 & 62.13 & 2 & 56.68 & 2 \\
80 & 50 & 75 & 1 & -0.9525 & -0.9525 & 102.11 & 1 & 22.94 & 1 \\
80 & 50 & 75 & 2 & -0.6360 & -0.7065 & 1,174.16 & 2 & 56.20 & 2 \\
80 & 50 & 75 & 3 & -0.8779 & -0.8779 & 240.58 & 2 & 49.64 & 2 \\
90 & 25 & 25 & 1 & -1.6797 & -1.6838 & 5.49 & 2 & 81.00 & 2 \\
90 & 25 & 25 & 2 & -1.2212 & -1.3070 & 6.52 & 2 & 55.72 & 2 \\
90 & 25 & 25 & 3 & -1.6920 & -1.6921 & 2.33 & 2 & 29.62 & 2 \\
90 & 25 & 50 & 1 & -1.0890 & -1.0890 & 99.33 & 2 & 55.78 & 2 \\
90 & 25 & 50 & 2 & -0.7685 & -0.7685 & 61.62 & 1 & 49.43 & 1 \\
90 & 25 & 50 & 3 & -1.0729 & -1.0807 & 66.62 & 2 & 65.67 & 2 \\
90 & 25 & 75 & 1 & -0.8472 & -0.8600 & 3,600.00 & 3 & 109.53 & 3 \\
90 & 25 & 75 & 2 & -0.7244 & -0.7244 & 1,426.76 & 1 & 48.41 & 1 \\
90 & 25 & 75 & 3 & -0.8361 & -0.8510 & 2,434.38 & 2 & 72.00 & 2 \\
90 & 50 & 25 & 1 & -1.6797 & -1.6838 & 5.16 & 2 & 74.51 & 2 \\
90 & 50 & 25 & 2 & -1.2943 & -1.3070 & 6.26 & 2 & 59.35 & 2 \\
90 & 50 & 25 & 3 & -1.6921 & -1.6921 & 1.25 & 1 & 14.62 & 1 \\
90 & 50 & 50 & 1 & -1.0511 & -1.0890 & 162.83 & 3 & 76.08 & 3 \\
90 & 50 & 50 & 2 & -0.7312 & -0.7685 & 130.59 & 2 & 89.46 & 2 \\
90 & 50 & 50 & 3 & -1.0690 & -1.0807 & 62.89 & 2 & 63.38 & 2 \\
90 & 50 & 75 & 1 & -0.8471 & -0.8598 & 2,413.71 & 2 & 69.97 & 2 \\
90 & 50 & 75 & 2 & -0.7223 & -0.7244 & 2,247.72 & 2 & 97.29 & 2 \\
90 & 50 & 75 & 3 & -0.8287 & -0.8510 & 3,600.00 & 1 & 73.69 & 1 \\
100 & 25 & 25 & 1 & -1.3988 & -1.4159 & 24.17 & 2 & 82.72 & 2 \\
100 & 25 & 25 & 2 & -1.1799 & -1.1849 & 21.10 & 2 & 87.04 & 2 \\
100 & 25 & 25 & 3 & -1.4414 & -1.4546 & 11.64 & 2 & 81.31 & 2 \\
100 & 25 & 50 & 1 & -0.9344 & -0.9503 & 2,976.25 & 3 & 396.85 & 3 \\
100 & 25 & 50 & 2 & -0.7645 & -0.8115 & 395.32 & 2 & 138.82 & 2 \\
100 & 25 & 50 & 3 & -1.0501 & -1.0604 & 163.69 & 2 & 77.10 & 2 \\
100 & 25 & 75 & 1 & -0.8683 & -0.8683 & 3,046.31 & 2 & 97.84 & 2 \\
100 & 25 & 75 & 2 & -0.6496 & -0.6672 & 3,600.00 & 1 & 486.76 & 3 \\
100 & 25 & 75 & 3 & -0.8128 & -0.8322 & 2,251.29 & 2 & 240.25 & 2 \\
100 & 50 & 25 & 1 & -1.4159 & -1.4159 & 10.82 & 1 & 43.00 & 1 \\
100 & 50 & 25 & 2 & -1.1609 & -1.1849 & 20.19 & 2 & 85.15 & 2 \\
100 & 50 & 25 & 3 & -1.4387 & -1.4546 & 12.52 & 2 & 89.28 & 2 \\
100 & 50 & 50 & 1 & -0.9344 & -0.9503 & 1,558.10 & 3 & 385.21 & 3 \\
100 & 50 & 50 & 2 & -0.8115 & -0.8115 & 158.51 & 1 & 76.04 & 1 \\
100 & 50 & 50 & 3 & -1.0245 & -1.0604 & 126.03 & 2 & 82.71 & 2 \\
100 & 50 & 75 & 1 & -0.8683 & -0.8683 & 1,962.57 & 2 & 101.57 & 2 \\
100 & 50 & 75 & 2 & -0.6526 & -0.6672 & 3,600.00 & 1 & 423.47 & 3 \\
100 & 50 & 75 & 3 & -0.8267 & -0.8322 & 2,177.76 & 2 & 185.07 & 2 \\
110 & 25 & 25 & 1 & -1.2951 & -1.3003 & 148.33 & 2 & 211.01 & 2 \\
110 & 25 & 25 & 2 & -1.1150 & -1.1412 & 46.44 & 2 & 115.67 & 2 \\
110 & 25 & 25 & 3 & -1.3923 & -1.4154 & 28.87 & 2 & 80.75 & 2 \\
110 & 25 & 50 & 1 & -0.8544 & -0.8781 & 3,537.69 & 2 & 377.36 & 2 \\
110 & 25 & 50 & 2 & -0.8046 & -0.8056 & 1,967.05 & 2 & 115.11 & 2 \\
110 & 25 & 50 & 3 & -0.9247 & -0.9392 & 2,381.22 & 3 & 581.00 & 3 \\
110 & 25 & 75 & 1 & -0.9423 & -0.9423 & 3,600.00 & 1 & 84.57 & 2 \\
110 & 25 & 75 & 2 & -0.7539 & -0.7599 & 3,600.00 & 1 & 382.18 & 2 \\
110 & 25 & 75 & 3 & -0.7741 & -0.7772 & 3,600.00 & 1 & 339.11 & 2 \\
110 & 50 & 25 & 1 & -1.2990 & -1.3003 & 1,516.03 & 2 & 151.64 & 2 \\
110 & 50 & 25 & 2 & -1.1402 & -1.1412 & 1,655.45 & 2 & 91.43 & 2 \\
110 & 50 & 25 & 3 & -1.3855 & -1.4154 & 40.99 & 2 & 72.87 & 2 \\
110 & 50 & 50 & 1 & -0.8765 & -0.8781 & 3,600.00 & 1 & 220.72 & 2 \\
110 & 50 & 50 & 2 & -0.7957 & -0.8056 & 3,600.00 & 1 & 92.10 & 2 \\
110 & 50 & 50 & 3 & -0.9319 & -0.9392 & 3,600.00 & 1 & 322.67 & 2 \\
110 & 50 & 75 & 1 & -0.9363 & -0.9423 & 3,600.00 & 1 & 111.70 & 3 \\
110 & 50 & 75 & 2 & -0.7539 & -0.7599 & 3,600.00 & 1 & 388.42 & 2 \\
110 & 50 & 75 & 3 & -0.7670 & -0.7772 & 3,600.00 & 1 & 368.62 & 2 \\
120 & 25 & 25 & 1 & -1.4269 & -1.4434 & 678.73 & 3 & 49.32 & 3 \\
120 & 25 & 25 & 2 & -1.1317 & -1.1398 & 2,009.52 & 3 & 156.35 & 3 \\
120 & 25 & 25 & 3 & -1.2517 & -1.2645 & 3,600.00 & 1 & 438.55 & 2 \\
120 & 25 & 50 & 1 & -0.9547 & -0.9547 & 3,600.00 & 1 & 59.06 & 1 \\
120 & 25 & 50 & 2 & -0.7075 & -0.7205 & 3,600.00 & 1 & 295.28 & 2 \\
120 & 25 & 50 & 3 & -0.8981 & -0.9082 & 3,600.00 & 1 & 480.57 & 2 \\
120 & 25 & 75 & 1 & -0.8806 & -0.8806 & 3,600.00 & 1 & 44.30 & 1 \\
120 & 25 & 75 & 2 & -0.6518 & -0.6535 & 3,600.00 & 1 & 911.78 & 3 \\
120 & 25 & 75 & 3 & -0.7976 & -0.8070 & 3,600.00 & 1 & 176.83 & 2 \\
120 & 50 & 25 & 1 & -1.4428 & -1.4434 & 35.67 & 2 & 113.19 & 2 \\
120 & 50 & 25 & 2 & -1.1344 & -1.1398 & 2,333.76 & 3 & 118.63 & 3 \\
120 & 50 & 25 & 3 & -1.2267 & -1.2645 & 2,719.24 & 2 & 302.04 & 2 \\
120 & 50 & 50 & 1 & -0.9544 & -0.9547 & 3,600.00 & 1 & 94.74 & 2 \\
120 & 50 & 50 & 2 & -0.7161 & -0.7205 & 3,600.00 & 1 & 352.88 & 2 \\
120 & 50 & 50 & 3 & -0.9065 & -0.9082 & 3,600.00 & 1 & 355.51 & 2 \\
120 & 50 & 75 & 1 & -0.8806 & -0.8806 & 3,600.00 & 1 & 34.54 & 1 \\
120 & 50 & 75 & 2 & -0.6456 & -0.6535 & 3,600.00 & 1 & 657.61 & 3 \\
120 & 50 & 75 & 3 & -0.7810 & -0.8070 & 3,600.00 & 1 & 166.77 & 2 \\
\end{longtable}
}

\begin{figure}[!ht]
    \centering
    \includegraphics[scale=0.6]{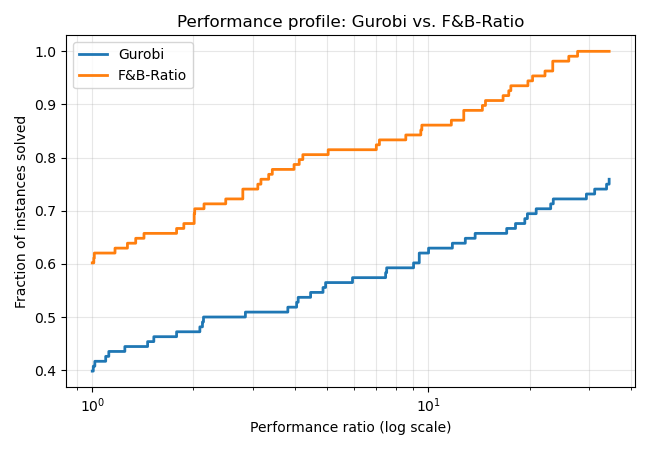}
    \caption{Performance profiles comparing Gurobi and \texttt{F\&B-Ratio} on the benchmark instances. A run is successful if its runtime is below 3600 seconds. The performance ratio is computed with respect to the best successful method on each instance and is shown on a logarithmic scale.}
    \label{fig:perf}
\end{figure}

\begin{figure}[!ht]
    \centering
    \includegraphics[scale=0.5]{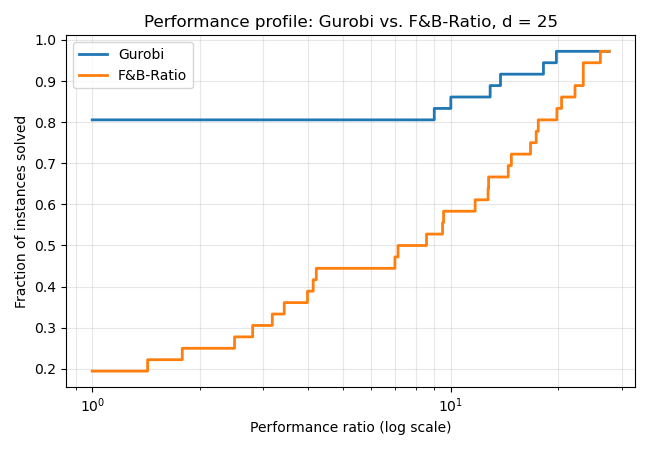}
    \includegraphics[scale=0.5]{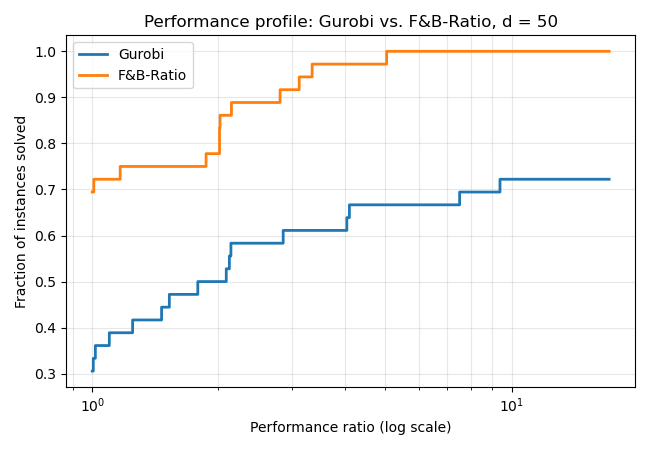}
    \includegraphics[scale=0.5]{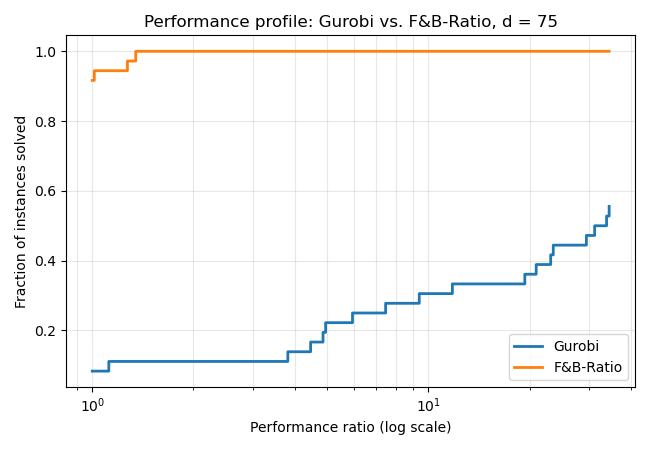}
    \caption{Performance profiles comparing Gurobi and \texttt{F\&B-Ratio} on the benchmark instances for densities $d=25$, $d=50$, and $d=75$.}
    \label{fig:perf_individual}
\end{figure}

The performance profiles in Figure~\ref{fig:perf} summarize the relative
efficiency of the two methods, while the density-wise profiles in
Figure~\ref{fig:perf_individual} provides a more detailed view of how this
behavior changes across instance classes. The horizontal axis reports the
performance ratio with respect to the best successful method on each instance,
using a logarithmic scale, while the vertical axis reports the fraction of
instances solved within a given ratio. Overall, the profiles confirm that Gurobi
is competitive on easier instances, where its overhead is smaller and it often
solves the subproblems very quickly. However, \texttt{F\&B-Ratio} dominates on
the harder instances: it solves a larger fraction of the test set within any
moderate performance ratio and is the only method that solves all instances
within the time limit. The density-wise profiles further show that this advantage
becomes more pronounced as the instances become denser.

{In Table~\ref{tab:init-detailed}, we compare the proposed heuristic
initialization with a random feasible initialization. The results clearly
show the benefit of the heuristic warm start. With the proposed
initialization, \texttt{F\&B-Ratio} requires 2 Dinkelbach iterations
on average, compared with 2.83 iterations when starting from a random
feasible point. Random initialization requires more iterations on 76
out of the 108 instances and the same number on the remaining 32, while
it never requires fewer iterations. The median number of iterations also
increases from 2 to 3, and the maximum from 3 to 4. These results confirm
that the high-quality initial solution produced by the hybrid heuristic
reduces the number of expensive global subproblems that must be solved.}

{
\footnotesize
\setlength{\tabcolsep}{4pt}
\begin{longtable}{rrrrrrrrrr}
\caption{{Detailed comparison of \texttt{F\&B-Ratio}
with the proposed heuristic initialization and with random initialization.
For each instance, the table reports the initial objective value, total
runtime in seconds, and number of Dinkelbach iterations. Both variants
converge to the same globally optimal objective value; therefore, the
latter is omitted.}}
\label{tab:init-detailed}\\

\toprule
\multicolumn{4}{c}{Instance}
& \multicolumn{3}{c}{Heuristic initialization}
& \multicolumn{3}{c}{Random initialization} \\
\cmidrule(lr){1-4}
\cmidrule(lr){5-7}
\cmidrule(lr){8-10}
$n$ & $p$ & $d$ & $s$
& Init. obj. & Time (s) & Iter.
& Init. obj. & Time (s) & Iter. \\
\midrule
\endfirsthead

\toprule
\multicolumn{4}{c}{Instance}
& \multicolumn{3}{c}{Heuristic initialization}
& \multicolumn{3}{c}{Random initialization} \\
\cmidrule(lr){1-4}
\cmidrule(lr){5-7}
\cmidrule(lr){8-10}
$n$ & $p$ & $d$ & $s$
& Init. obj. & Time (s) & Iter.
& Init. obj. & Time (s) & Iter. \\
\midrule
\endhead

\midrule
\multicolumn{10}{r}{Continued on next page}\\
\endfoot

\bottomrule
\endlastfoot

70 & 25 & 25 & 1 & -2.1005 & 20.22 & 2 & -0.0384 & 51.55 & 3 \\
70 & 25 & 25 & 2 & -1.3065 & 49.96 & 2 & -0.1033 & 106.29 & 3 \\
70 & 25 & 25 & 3 & -1.9588 & 10.31 & 1 & -0.2005 & 36.26 & 3 \\
70 & 25 & 50 & 1 & -1.1667 & 76.17 & 3 & -0.0201 & 75.68 & 3 \\
70 & 25 & 50 & 2 & -0.9003 & 39.80 & 2 & 0.2119 & 44.82 & 2 \\
70 & 25 & 50 & 3 & -1.1748 & 35.23 & 2 & -0.1358 & 58.24 & 3 \\
70 & 25 & 75 & 1 & -0.9756 & 110.11 & 2 & -0.1597 & 145.09 & 3 \\
70 & 25 & 75 & 2 & -0.6321 & 96.37 & 2 & 0.1097 & 118.68 & 2 \\
70 & 25 & 75 & 3 & -1.0203 & 53.80 & 2 & 0.0126 & 40.18 & 2 \\
70 & 50 & 25 & 1 & -2.1564 & 19.06 & 2 & -0.0906 & 29.11 & 3 \\
70 & 50 & 25 & 2 & -1.2877 & 51.06 & 2 & -0.2056 & 74.16 & 3 \\
70 & 50 & 25 & 3 & -1.9588 & 11.00 & 1 & 0.2631 & 33.14 & 3 \\
70 & 50 & 50 & 1 & -1.1712 & 81.11 & 2 & -0.0212 & 71.28 & 3 \\
70 & 50 & 50 & 2 & -0.8693 & 38.39 & 2 & 0.1259 & 39.65 & 2 \\
70 & 50 & 50 & 3 & -1.1905 & 35.46 & 2 & -0.0171 & 56.55 & 3 \\
70 & 50 & 75 & 1 & -0.9765 & 77.32 & 2 & -0.0994 & 105.62 & 3 \\
70 & 50 & 75 & 2 & -0.6851 & 98.15 & 2 & 0.0494 & 189.27 & 2 \\
70 & 50 & 75 & 3 & -1.0203 & 46.80 & 2 & -0.1528 & 85.55 & 2 \\

\midrule

80 & 25 & 25 & 1 & -1.8984 & 25.88 & 3 & -0.1285 & 102.89 & 3 \\
80 & 25 & 25 & 2 & -1.4282 & 55.06 & 2 & -0.2361 & 97.94 & 3 \\
80 & 25 & 25 & 3 & -1.8923 & 12.89 & 1 & -0.0166 & 48.37 & 3 \\
80 & 25 & 50 & 1 & -1.1612 & 35.07 & 1 & -0.1305 & 97.59 & 3 \\
80 & 25 & 50 & 2 & -0.8626 & 24.02 & 1 & 0.2221 & 80.31 & 3 \\
80 & 25 & 50 & 3 & -1.1478 & 56.89 & 2 & -0.1961 & 96.50 & 3 \\
80 & 25 & 75 & 1 & -0.9207 & 48.84 & 2 & -0.0176 & 100.10 & 3 \\
80 & 25 & 75 & 2 & -0.6802 & 54.32 & 2 & 0.0264 & 113.04 & 3 \\
80 & 25 & 75 & 3 & -0.8620 & 76.47 & 3 & -0.0439 & 91.55 & 4 \\
80 & 50 & 25 & 1 & -1.9175 & 12.60 & 1 & -0.2442 & 47.82 & 3 \\
80 & 50 & 25 & 2 & -1.4151 & 56.30 & 2 & 0.0671 & 81.08 & 3 \\
80 & 50 & 25 & 3 & -1.8901 & 25.06 & 2 & -0.0320 & 40.80 & 3 \\
80 & 50 & 50 & 1 & -1.1547 & 73.04 & 2 & -0.1793 & 93.23 & 3 \\
80 & 50 & 50 & 2 & -0.8110 & 47.98 & 2 & 0.1906 & 75.45 & 3 \\
80 & 50 & 50 & 3 & -1.1277 & 56.68 & 2 & 0.0033 & 82.13 & 3 \\
80 & 50 & 75 & 1 & -0.9525 & 22.94 & 1 & -0.0826 & 73.71 & 3 \\
80 & 50 & 75 & 2 & -0.6360 & 56.20 & 2 & 0.0099 & 95.64 & 2 \\
80 & 50 & 75 & 3 & -0.8779 & 49.64 & 2 & 0.0192 & 101.02 & 4 \\

\midrule

90 & 25 & 25 & 1 & -1.6797 & 81.00 & 2 & -0.1148 & 109.24 & 3 \\
90 & 25 & 25 & 2 & -1.2212 & 55.72 & 2 & 0.0211 & 88.62 & 3 \\
90 & 25 & 25 & 3 & -1.6920 & 29.62 & 2 & -0.1008 & 41.45 & 3 \\
90 & 25 & 50 & 1 & -1.0890 & 55.78 & 2 & -0.0268 & 73.54 & 3 \\
90 & 25 & 50 & 2 & -0.7685 & 49.43 & 1 & 0.2251 & 176.82 & 3 \\
90 & 25 & 50 & 3 & -1.0729 & 65.67 & 2 & 0.0921 & 195.91 & 3 \\
90 & 25 & 75 & 1 & -0.8472 & 109.53 & 3 & -0.1009 & 129.96 & 3 \\
90 & 25 & 75 & 2 & -0.7244 & 48.41 & 1 & -0.0124 & 118.68 & 2 \\
90 & 25 & 75 & 3 & -0.8361 & 72.00 & 2 & 0.0367 & 104.57 & 3 \\
90 & 50 & 25 & 1 & -1.6797 & 74.51 & 2 & -0.1780 & 111.72 & 3 \\
90 & 50 & 25 & 2 & -1.2943 & 59.35 & 2 & -0.1582 & 86.88 & 3 \\
90 & 50 & 25 & 3 & -1.6921 & 14.62 & 1 & 0.0673 & 43.79 & 3 \\
90 & 50 & 50 & 1 & -1.0511 & 76.08 & 3 & -0.0003 & 71.80 & 3 \\
90 & 50 & 50 & 2 & -0.7312 & 89.46 & 2 & 0.0313 & 128.69 & 3 \\
90 & 50 & 50 & 3 & -1.0690 & 63.38 & 2 & 0.0850 & 98.22 & 3 \\
90 & 50 & 75 & 1 & -0.8471 & 69.97 & 2 & -0.1500 & 110.11 & 3 \\
90 & 50 & 75 & 2 & -0.7223 & 97.29 & 2 & 0.0264 & 122.78 & 2 \\
90 & 50 & 75 & 3 & -0.8287 & 73.69 & 1 & 0.0260 & 177.94 & 3 \\

\midrule

100 & 25 & 25 & 1 & -1.3988 & 82.72 & 2 & -0.1223 & 127.98 & 2 \\
100 & 25 & 25 & 2 & -1.1799 & 87.04 & 2 & 0.0426 & 93.72 & 2 \\
100 & 25 & 25 & 3 & -1.4414 & 81.31 & 2 & 0.0812 & 84.47 & 2 \\
100 & 25 & 50 & 1 & -0.9344 & 396.85 & 3 & -0.0196 & 677.09 & 4 \\
100 & 25 & 50 & 2 & -0.7645 & 138.82 & 2 & 0.1411 & 201.83 & 3 \\
100 & 25 & 50 & 3 & -1.0501 & 77.10 & 2 & -0.1090 & 90.29 & 2 \\
100 & 25 & 75 & 1 & -0.8683 & 97.84 & 2 & -0.0649 & 184.83 & 3 \\
100 & 25 & 75 & 2 & -0.6496 & 486.76 & 3 & 0.0637 & 508.00 & 3 \\
100 & 25 & 75 & 3 & -0.8128 & 240.25 & 2 & 0.0305 & 201.25 & 3 \\
100 & 50 & 25 & 1 & -1.4159 & 43.00 & 1 & -0.1351 & 130.06 & 2 \\
100 & 50 & 25 & 2 & -1.1609 & 85.15 & 2 & -0.1661 & 89.16 & 2 \\
100 & 50 & 25 & 3 & -1.4387 & 89.28 & 2 & -0.0348 & 107.36 & 2 \\
100 & 50 & 50 & 1 & -0.9344 & 385.21 & 3 & -0.1045 & 768.57 & 4 \\
100 & 50 & 50 & 2 & -0.8115 & 76.04 & 1 & 0.0568 & 233.52 & 3 \\
100 & 50 & 50 & 3 & -1.0245 & 82.71 & 2 & -0.0754 & 96.02 & 2 \\
100 & 50 & 75 & 1 & -0.8683 & 101.57 & 2 & -0.1366 & 202.18 & 3 \\
100 & 50 & 75 & 2 & -0.6526 & 423.47 & 3 & 0.0122 & 613.30 & 3 \\
100 & 50 & 75 & 3 & -0.8267 & 185.07 & 2 & 0.0012 & 247.97 & 3 \\

\midrule

110 & 25 & 25 & 1 & -1.2951 & 211.01 & 2 & 0.0205 & 511.33 & 3 \\
110 & 25 & 25 & 2 & -1.1150 & 115.67 & 2 & -0.0538 & 157.82 & 3 \\
110 & 25 & 25 & 3 & -1.3923 & 80.75 & 2 & -0.1261 & 143.45 & 3 \\
110 & 25 & 50 & 1 & -0.8544 & 377.36 & 2 & 0.0147 & 620.55 & 3 \\
110 & 25 & 50 & 2 & -0.8046 & 115.11 & 2 & 0.1140 & 101.04 & 2 \\
110 & 25 & 50 & 3 & -0.9247 & 581.00 & 3 & -0.1563 & 602.19 & 3 \\
110 & 25 & 75 & 1 & -0.9423 & 84.57 & 2 & -0.0898 & 188.28 & 4 \\
110 & 25 & 75 & 2 & -0.7539 & 382.18 & 2 & -0.0201 & 435.69 & 2 \\
110 & 25 & 75 & 3 & -0.7741 & 339.11 & 2 & -0.0655 & 608.92 & 3 \\
110 & 50 & 25 & 1 & -1.2990 & 151.64 & 2 & -0.0861 & 450.04 & 3 \\
110 & 50 & 25 & 2 & -1.1402 & 91.43 & 2 & 0.0626 & 200.99 & 3 \\
110 & 50 & 25 & 3 & -1.3855 & 72.87 & 2 & 0.0323 & 264.13 & 3 \\
110 & 50 & 50 & 1 & -0.8765 & 220.72 & 2 & -0.0453 & 511.52 & 3 \\
110 & 50 & 50 & 2 & -0.7957 & 92.10 & 2 & -0.0341 & 192.54 & 2 \\
110 & 50 & 50 & 3 & -0.9319 & 322.67 & 2 & -0.0068 & 295.30 & 2 \\
110 & 50 & 75 & 1 & -0.9363 & 111.70 & 3 & -0.2018 & 323.25 & 4 \\
110 & 50 & 75 & 2 & -0.7539 & 388.42 & 2 & -0.0504 & 332.10 & 2 \\
110 & 50 & 75 & 3 & -0.7670 & 368.62 & 2 & -0.0285 & 581.10 & 3 \\

\midrule

120 & 25 & 25 & 1 & -1.4269 & 49.32 & 3 & -0.2439 & 128.38 & 3 \\
120 & 25 & 25 & 2 & -1.1317 & 156.35 & 3 & -0.1410 & 214.72 & 4 \\
120 & 25 & 25 & 3 & -1.2517 & 438.55 & 2 & -0.1286 & 716.38 & 3 \\
120 & 25 & 50 & 1 & -0.9547 & 59.06 & 1 & -0.0132 & 387.04 & 3 \\
120 & 25 & 50 & 2 & -0.7075 & 295.28 & 2 & 0.0110 & 330.69 & 2 \\
120 & 25 & 50 & 3 & -0.8981 & 480.57 & 2 & -0.0122 & 798.94 & 3 \\
120 & 25 & 75 & 1 & -0.8806 & 44.30 & 1 & -0.0513 & 108.62 & 2 \\
120 & 25 & 75 & 2 & -0.6518 & 911.78 & 3 & -0.0102 & 1,611.94 & 4 \\
120 & 25 & 75 & 3 & -0.7976 & 176.83 & 2 & -0.0378 & 196.29 & 2 \\
120 & 50 & 25 & 1 & -1.4428 & 113.19 & 2 & -0.0852 & 109.30 & 3 \\
120 & 50 & 25 & 2 & -1.1344 & 118.63 & 3 & 0.0566 & 266.35 & 4 \\
120 & 50 & 25 & 3 & -1.2267 & 302.04 & 2 & 0.0443 & 565.14 & 3 \\
120 & 50 & 50 & 1 & -0.9544 & 94.74 & 2 & -0.0659 & 373.77 & 3 \\
120 & 50 & 50 & 2 & -0.7161 & 352.88 & 2 & 0.0768 & 379.62 & 2 \\
120 & 50 & 50 & 3 & -0.9065 & 355.51 & 2 & 0.0098 & 797.91 & 3 \\
120 & 50 & 75 & 1 & -0.8806 & 34.54 & 1 & -0.0530 & 103.49 & 2 \\
120 & 50 & 75 & 2 & -0.6456 & 657.61 & 3 & 0.0005 & 1,219.99 & 4 \\
120 & 50 & 75 & 3 & -0.7810 & 166.77 & 2 & -0.1158 & 171.42 & 2 \\

\end{longtable}
}

{To further assess the scalability of \texttt{F\&B-Ratio}, we generated an additional set of 18 larger instances with $n=160$, using the same instance-generation procedure described in Section \ref{sec:inst}. The same time limit of 3,600 seconds was imposed. The results are reported in Table \ref{tab:large_instances}. In contrast to the results reported in Table \ref{tab:gurobi-fnb-full}, \texttt{F\&B-Ratio} does not solve all these larger instances within the time limit: 13 out of 18 instances are solved, while the remaining five reach the time limit. All instances with $d=25$ are solved, whereas two instances with $d=50$ reach the time limit. The most challenging setting is $p=50$ and $d=75$, for which all three instances reach the time limit. These results show that $n=160$ is the maximum size our method can handle in a reasonable amount of time.}

\begin{table}[!ht]
\centering
\caption{{Results of \texttt{F\&B-Ratio} on the larger
instances with $n=160$. For each instance, the table reports the initial
objective value, the objective value returned by \texttt{F\&B-Ratio} if
it proves optimality, a dash otherwise, the total runtime in seconds,
and the number of Dinkelbach iterations. A runtime of 3,600 seconds
indicates that the time limit was reached.}}
\label{tab:large_instances}

\small
\renewcommand{\arraystretch}{1.18}
\setlength{\tabcolsep}{3.5pt}

\begin{tabular}{rrrrrrrr}
\hline
\multicolumn{4}{c}{Instance} &
\multicolumn{4}{c}{\texttt{F\&B-Ratio}} \\[2pt]
$n$ & $p$ & $d$ & $s$ &
Init. obj. & Obj. & Time (s) & Iter. \\[2pt]
\hline
160 & 25 & 25 & 1 & -1.1467 & -1.1473 &  616.39 & 2 \\
160 & 25 & 25 & 2 & -0.9462 & -0.9755 &  494.93 & 2 \\
160 & 25 & 25 & 3 & -1.1405 & -1.1503 &  740.66 & 2 \\
160 & 25 & 50 & 1 & -0.8303 & -0.8502 &  154.90 & 2 \\
160 & 25 & 50 & 2 & -0.6596 & --      & 3600.00 & 2 \\
160 & 25 & 50 & 3 & -0.8317 & -0.8581 &  584.34 & 2 \\
160 & 25 & 75 & 1 & -0.7103 & -0.7184 & 1537.31 & 2 \\
160 & 25 & 75 & 2 & -0.6301 & -0.6418 & 1095.63 & 3 \\
160 & 25 & 75 & 3 & -0.6704 & -0.6813 & 2147.63 & 2 \\
\hline
160 & 50 & 25 & 1 & -1.1467 & -1.1473 &  704.83 & 2 \\
160 & 50 & 25 & 2 & -0.9688 & -0.9755 &  497.67 & 2 \\
160 & 50 & 25 & 3 & -1.1350 & -1.1503 &  835.19 & 2 \\
160 & 50 & 50 & 1 & -0.8477 & -0.8502 &  166.63 & 2 \\
160 & 50 & 50 & 2 & -0.6532 & --      & 3600.00 & 2 \\
160 & 50 & 50 & 3 & -0.8370 & -0.8581 &  507.01 & 2 \\
160 & 50 & 75 & 1 & -0.7134 & --      & 3600.00 & 2 \\
160 & 50 & 75 & 2 & -0.6402 & --      & 3600.00 & 1 \\
160 & 50 & 75 & 3 & -0.6717 & --      & 3600.00 & 1 \\
\hline
\end{tabular}
\end{table}

\section{Conclusions}\label{sec:concl}
In this paper, we studied mixed-integer bound-constrained quadratic fractional
optimization problems, where the objective is given by the ratio of two
quadratic functions and the feasible set contains both binary and continuous
box-constrained variables. This class of problems is challenging because it
combines the nonconvexity of fractional quadratic objectives with mixed-integer
constraints.

We proposed an exact solution approach based on Dinkelbach's parametric
reformulation. The original fractional problem is solved by iteratively solving
a sequence of parametric nonconvex MIQP subproblems to global optimality. To
handle these subproblems efficiently, we exploit a tailored solver that
reformulates them as box-constrained quadratic programs. In addition, we introduced a hybrid heuristic based on continuous local search and Variable Neighborhood Search to compute
high-quality initial feasible solutions for the Dinkelbach scheme.
The computational results confirm the effectiveness of the proposed method, which outperforms the general-purpose MIQP solver on the generated benchmark set, particularly on the most challenging instances. Moreover, the heuristic provides high-quality initial solutions, leading to a small number of Dinkelbach iterations and improving the efficiency of the overall exact scheme.

Future work will focus on extending the computational study to larger-scale
instances and additional classes of fractional mixed-integer quadratic problems.

\section*{Acknowledgments} We thank the anonymous referee for the careful reading of the paper and for the suggestions that greatly helped us to improve the paper.

%
%
%

\end{document}